\documentclass[12pt,reqno]{amsart}
\usepackage{amsmath}
\usepackage{amssymb}
\usepackage{latexsym}

\newcommand{\Zint}{\mathbb {Z}}    
     
\newcommand{\Rea}{\mathbb {R}}      

\newcommand{\halmos}{\rule{5pt}{5pt}}

\numberwithin{equation}{section}

\newtheorem{prop}{\bf Proposition}[section]
\newtheorem{thm}[prop]{\bf Theorem}

\newtheorem{exa}{\bf Example}

\begin{document}
\baselineskip 15pt

\title[On eigenvalues of Lam\'e operator]
{On eigenvalues of Lam\'e operator}
\author{Kouichi Takemura}
\address{Department of Mathematical Sciences, Yokohama City University, 22-2 Seto, Kanazawa-ku, Yokohama 236-0027, Japan.}
\email{takemura@yokohama-cu.ac.jp}

\begin{abstract}
We introduce two integral representations of monodromy on Lam\'e equation. By applying them, we obtain results on hyperelliptic-to-elliptic reduction integral formulae, finite-gap potential and eigenvalues of Lam\'e operator.
\end{abstract}

\maketitle

\section{Introduction}

The Lam\'e operator is a differential operator defined by
\begin{equation}
H = -\frac{d ^2}{d x^2} +n(n+1) \wp (x) . \label{LameH}
\end{equation}
where $\wp (x)$ is the Weierstrass $\wp$-function and $n$ is a constant. In this manuscript we assume that basic periods of the function $\wp (x)$ are $(1,\tau)$ and $n$ is a positive integer.

We are going to consider eigenvalues of the Lam\'e operator $H$ with boundary condition. 
In other words, we specify the vector space where the operator $H$ acts and consider the spectral problem on the space. 

Let $f(x)$ be an eigenfunction of the operator $H$ with eigenvalue $E$, i.e.
\begin{equation}
(H-E) f(x) = \left( -\frac{d ^2}{d x^2} +n(n+1) \wp (x) -E \right) f(x) =0.
\label{eqn:lame}
\end{equation}
Then Eq.(\ref{eqn:lame}) is called the Lam\'e equation (or the Lam\'e's differential equation).
Although the differential equation (\ref{eqn:lame}) is periodic in $x$ with period $1$, the solution may not be periodic. We also note that the eigenfunction may not be square-integrable on the interval $(0,1)$ because the potential  $n(n+1)\wp (x)$ has poles on $\Zint $. In this manuscript we consider the condition for the eigenvalue $E$ such that Eq.(\ref{eqn:lame}) has a singly-periodic eigenfunction (i.e. $f(x+1)=\pm f(x)$), a doubly-periodic eigenfunction (i.e. $f(x+1)=\pm f(x)$ and $f(x+\tau )=\pm f(x)$) or a square-integrable eigenfunction (i.e. $\int _0 ^1 |f(x)|^2 dx < + \infty $). These conditions correspond to the boundary condition and the boundary condition is closely related with the monodromy. For the case $n \in \Zint $, the monodromy of solutions to Eq.(\ref{eqn:lame}) is calculated and it has two expressions. One is expressed by hyperelliptic integral and the other is based on Hermite-Krichever Ansatz. Note that hyperelliptic-to-elliptic reduction integral formulae are obtained by comparing two expressions.

A crucial fact on the Lam\'e equation is that the Lam\'e operator is an example of the finite-gap potential which was established by Ince \cite{I}. By applying the formula of the monodromy in terms of hyperelliptic integral, we can show results on finite-gap potential. Moreover relationship between the boundary condition and the finite-gap potential is clarified.

The eigenvalue $E$ of the  Lam\'e operator with boundary condition depends on the period $\tau$.
In particular we can numerically investigate branching of the eigenvalues as a complex-analytic function in $\tau $ for the case $n=1$, and it is compatible with the convergence radius of the eigenvalue expanded as a power series in $p=\exp (\pi \sqrt{-1} \tau )$.

Here we comment on relationship with quantum integrable system.
The elliptic Calogero-Moser-Sutherland model (or elliptic Olshanetsky-Perelomov model \cite{OP}) of type $A_N$ is a quantum many body system whose Hamiltonian is given as follows,
\begin{equation}
 H _{A_N}=- \frac{1}{2} \sum_{i=1}^{N+1} 
\frac{\partial ^{2}}{\partial x_{i}^{2}}
+ n (n+1)
\sum_{1 \leq i<j \leq N+1}
\wp ( x_{i}-x_{j}).
\end{equation}
This model is known to be integrable, i.e. there exists $N+1$-algebraically independent commuting operators which commute with the Hamiltonian $H_{A_N}$.
For the case $N=1$, we recover the Lam\'e operator by setting $x=x_1-x_2$ and $y=x_1+x_2$. Therefore the elliptic Calogero-Moser-Sutherland model is a generalization of the Lam\'e operator. 
In contract to the case of the trigonometric Calogero-Moser-Sutherland model, spectral problem of the elliptic Calogero-Moser-Sutherland model is not studied very much.
We hope that the spectral problem of the elliptic Calogero-Moser-Sutherland model is clarified by developing knowledge on the Lam\'e equation in future. The paper by Chalykh, Etingof and Oblomkov \cite{CEO} might suggest us an approach to this problem.

This manuscript is organized as follows.
In section \ref{sec:othereq}, we describe relationship among Lam\'e equation, Heun's differential equation \cite{Ron} and hypergeometric polynomial. 
In section \ref{sec:prod}, we introduce a doubly-periodic function that satisfies a differential equation of order three, and obtain an integral representation of solutions to the Lam\'e equation and a formula of global monodromy in terms of hyperelliptic integral.
In section \ref{sec:HKA}, we explain results on Bethe Ansatz and Hermite-Krichever Ansatz. As an application of Hermite-Krichever Ansatz, we get a formula of global monodromy in terms of elliptic integral.
In section \ref{sec:intf}, we obtain hyperelliptic-to-elliptic reduction integral formulae by comparing two expressions of monodromy.
In section \ref{sec:bound}, we describe results on the finite-gap potential.
In section \ref{sec:n1}, we consider analytic continuation of eigenvalues in the variable $\tau$ for the case $n=1$.

\section{Relationship with other equations} \label{sec:othereq}

It is known that Lam\'e equation is a special case of Heun equation.
The Heun equation (or the Heun's differential equation) is a differential equation given by
\begin{equation}
\left( \! \left(\frac{d}{dw}\right) ^2 \! + \left( \frac{\gamma}{w}+\frac{\delta }{w-1}+\frac{\epsilon}{w-t}\right) \frac{d}{dw} +\frac{\alpha \beta w -q}{w(w-1)(w-t)} \right)\tilde{f}(w)=0
\label{Heun}
\end{equation}
with the condition 
\begin{equation}
\gamma +\delta +\epsilon =\alpha +\beta +1.
\label{Heuncond}
\end{equation}
Note that Eq.(\ref{Heun}) has four singularities $\{ 0,1,t,\infty \}$ on the Riemann sphere and the Heun equation is the standard canonical form of a Fuchsian equation with four singularities.
It is well known that the Fuchsian equation with three singularities is the hypergeometric differential equation.

The Heun equation admits an expression by elliptic function. Let $H_H$ be the Hamiltonian of $BC_1$ Inozemtsev model which is given as
\begin{equation}
H_H= -\frac{d^2}{dx^2} + \sum_{i=0}^3 l_i(l_i+1)\wp (x+\omega_i),
\label{Ino}
\end{equation}
where $\omega _0=0$, $\omega_1=1/2$, $\omega_2=-(\tau+1)/2$ and $\omega_3=\tau/2$ are half-periods, and $l_i$ $(i=0,1,2,3)$ are coupling constants. This model is a one-particle version of the $BC_N$ Inozemtsev system \cite{Ino}, which is known to be the universal quantum integrable system with $B_N$ symmetry \cite{Ino,OOS}. Let $f(x)$ be an eigenfuction of the operator $H_H$ with eigenvalue $E$, i.e.
\begin{equation}
H_H f(x) = E f(x).
\label{InoEF2}
\end{equation}
Under the transformation 
\begin{equation}
w=\frac{e_1-e_3}{\wp (x)-e_3}, \quad \quad w^{\frac{l_0+1}{2}}(w-1)^{\frac{l_1+1}{2}}(w-t)^{\frac{l_2+1}{2}} \tilde{f}(w) = f(x),
\label{wxtrans}
\end{equation}
Eq.(\ref{InoEF2}) is transformed to the Heun equation (\ref{Heun}). In this sense, Eq.(\ref{InoEF2}) is an elliptic representation of the Heun equation. The coupling constants $l_0, \; l_1, \; l_2, \; l_3$ correspond to the exponents $\alpha , \dots , \epsilon$, the elliptic modurus $\tau$ corresponds to the singular point $t$ and the eigenvalue $E$ corresponds to the accessory parameter $q$. For details see \cite{Tak2}. If $l_0 \neq 0$ and $l_1=l_2=l_3=0$, then Eq.(\ref{InoEF2}) represents the Lam\'e equation. Thus the Lam\'e equation is a special case of the Heun equation.

We observe a relationship with hypergeometric polynomial.
More precisely, hypergeomertic differential equation is obtained from the Lam\'e equation (or the Heun equation) by trigonometric limit.
Set 
\begin{align}
& H_T=  -\frac{d^2}{dx^2} + n(n+1)\frac{\pi ^2}{\sin ^2\pi x} .\label{triIno} 
\end{align}
Then $H \rightarrow H_T -\frac{\pi^2}{3} n (n +1)$ by the trigonometric limit $p = \exp ( \pi \sqrt{-1} \tau )\rightarrow 0$.
The equation $(H_T -E)f(x)=0$ is transformed to the hypergeometric differential equation. Set
\begin{equation}
\tilde{\Phi}(x)=(\sin \pi x)^{n+1}, \; \; \; v_m=  \tilde{c}_m C^{n+1}_m (\cos \pi x) \tilde{\Phi}(x) ,\; \; \; (m \in \Zint_{\geq 0}),
\label{Gegenpol}
\end{equation}
where the function $C^{\nu}_m (z)=\frac{\Gamma (m+2\nu)}{m! \Gamma(2\nu)}\, _2 \! F_1(-m,m+2\nu ; \nu+\frac{1}{2} ;\frac{1-z}{2})$ is the Gegenbauer polynomial of degree $m$ and $\tilde{c} _m=\sqrt{\frac{2^{2n+1}(m+n+1)m! \Gamma (n+1)^2}{\Gamma (m+2n+2)}}$.
Then we have
\begin{align}
H_T v_m=\pi^2(m+n+1)^2 v_m,
\end{align}
and the set $\{ v_m \} _{m \in \Zint _{\geq 0}}$ forms a complete orthonormal system on a certain Hilbert space. Hence the spectral problem for the operator $H_T$ is clarified substantially.

To investigate the spectral problem for the Lam\'e operator $H$, we can apply a method of perturbation from the operator $H_T$ and we have an algorithm for obtaining eigenvalues and eigenfunctions of the operator $H$ as formal power series in $p$. It is shown in \cite{Tak,KT,Tak2} that, if $|p|$ is sufficiently small, then the formal power series of eigenvalues and eigenfunctions converge.
Relationship between the convergence radius and the monodromy will be mentioned later.

\section{Monodromy and hyperelliptic integral} \label{sec:prod}

In this section we review results on an integral representation of solutions and a monodromy formula which is expressed by hyperelliptic integral. For this purpose we introduce a doubly-periodic function which plays important roles.
Let $h(x)$ be the product of any pair of the solutions to Eq.(\ref{eqn:lame}). Then the function $h(x)$ satisfies the following third-order differential equation:
\begin{align}
& \left( \frac{d^3}{dx^3}-4\left( n(n+1)\wp (x)-E\right)\frac{d}{dx}-2 n(n+1)\wp '(x) \right) h (x)=0.
\label{prodDE} 
\end{align}

It is known that Eq.(\ref{prodDE}) has a nonzero doubly-periodic solution for all $E$.
\begin{prop} \cite[Proposition 3.5]{Tak1} \label{prop:prod}
If $n\in \Zint _{\geq 1}$, then Eq.(\ref{prodDE}) has a nonzero even doubly-periodic solution $\Xi (x,E)$, which has the expansion
\begin{equation}
\Xi (x,E)=c_0(E)+\sum_{j=0}^{n-1} b_j (E)\wp (x)^{n-j},
\label{Fx}
\end{equation}
where the coefficients $c_0(E)$ and $b_j(E)$ are polynomials in $E$, they do not have common divisors and the polynomial $c_0(E)$ is monic.
\end{prop}
For the case of Lam\'e equation, we have $\deg_E c_0(E)=n$ and $\deg _E b_j(E)=j$.

We can derive an integral formula for the solution to Eq.(\ref{eqn:lame}) by using the doubly-periodic function $\Xi(x,E)$. Set
\begin{align}
 & Q(E)=  \Xi (x,E)^2\left( E- n(n+1)\wp (x)\right) +\frac{1}{2}\Xi (x,E)\frac{d^2\Xi (x,E)}{dx^2}-\frac{1}{4}\left(\frac{d\Xi (x,E)}{dx} \right)^2. \label{const}
\end{align}
Then the r.h.s. is independent of $x$, and $Q(E)$ is a monic polynomial in $E$ of degree $2n+1$. The following proposition on the integral representation of solutions is obtained in \cite{Tak1}:
\begin{prop} \cite[Proposition 3.7]{Tak1}
Let $\Xi (x,E)$ be the doubly-periodic function defined in Proposition \ref{prop:prod} and $Q(E)$ be the monic polynomial defined in Eq.(\ref{const}).
Then the function 
\begin{equation}
\Lambda (x,E)=\sqrt{\Xi (x,E)}\exp \int \frac{\sqrt{-Q(E)}dx}{\Xi (x,E)}
\label{eqn:Lam}
\end{equation}
is a solution to the differential equation (\ref{eqn:lame}).
\end{prop}
Note that Eq.(\ref{eqn:Lam}) is written in the book of Whittaker and Watson \cite{WW} for the case of the Lam\'e equation.

\begin{exa}
(i) The case $n=1$.
\begin{align}
& \Xi (x,E)= \wp (x) +E , \quad Q(E)=(E+e_1)(E+e_2)(E+e_3).
\end{align}
(ii) The case $n=2$.
\begin{align}
& \textstyle \Xi (x,E)= 9\wp (x)^2 +3E \wp (x) +E^2 -\frac{9}{4} g_2, \quad Q(E)=(E^2 -3g_2)\prod _{i=1}^3 (E-3e_i),
\end{align}
where $g_2= -4(e_1e_2+e_2e_3+e_3e_1)$.\\
(iii) The case $n=3$.
\begin{align}
& \textstyle \Xi (x,E)= 225\wp(x)^3+45E\wp(x)^2+6(E^2-\frac{75}{8}g_2)\wp(x)+E^3-15g_2 E-\frac{225}{4}g_3, \\
& \textstyle Q(E)= E\prod _{i=1}^3 (E^2+6e_i E+45e_i ^2-15g_2),
\end{align}
where $g_3=4e_1e_2e_3$.
\end{exa}

Now we consider the monodromy. It follows from Eq.(\ref{eqn:Lam}) that $\Lambda (x +1,E)$ (resp. $\Lambda (x +\tau ,E)$) is expressed as $m_1  \Lambda (x,E)$ (resp. $m_{\tau } \Lambda (x ,E)$) for some constant $m_1$ (resp. $m_{\tau }$). We determine the constants $m_1$ and $m_{\tau }$.
Rewrite the function $\Xi (x,E)$ as follows:
\begin{equation}
\Xi (x,E)=c(E)+\sum_{j=0}^{n-1 } a_j (E)\left( \frac{d}{dx} \right) ^{2j} \wp (x).
\label{FFx}
\end{equation}
\begin{prop} \cite{Tak3,Tak4} \ \label{thm:conj3} 
Assume $n \in \Zint_{\geq 1}$. Let $E'$ be the eigenvalue satisfying $Q(E')=0$. Then there exists $q_1, q_3 \in \{0,1\}$ such that $\Lambda (x+2\omega _k,E')=(-1)^{q_k} \Lambda (x,E')$, and we have
\begin{align} 
& \Lambda (x+1,E)=(-1)^{q_1} \Lambda (x,E) \exp \left( -\frac{1}{2} \int_{E'}^{E}\frac{c(\tilde{E})  -2\eta _1 a_0 (\tilde{E}) }{\sqrt{-Q(\tilde{E})}} d\tilde{E}\right) ,
\label{hypellint} \\
& \Lambda (x+\tau,E)=(-1)^{q_3} \Lambda (x,E) \exp \left( -\frac{1}{2} \int_{E'}^{E}\frac{\tau c(\tilde{E})  -2\eta _3 a_0 (\tilde{E}) }{\sqrt{-Q(\tilde{E})}} d\tilde{E}\right) ,
\label{hypellinttau}
\end{align}
where $\eta _1 =\zeta (1/2)$, $\eta _3 =\zeta (\tau/2)$ and $\zeta (x)$ is the Weierstrass zeta function.
\end{prop}
\begin{exa}
(i) The case $n=1$. Set $E'=-e_1$. Then we have $q_1 =0$, $q_3=1$. The polynomials $c(E)$ and $a_0 (E) $ are calculated as
\begin{align}
& c (E)= E , \quad a_0 (E) = 1.
\end{align}
(ii) The case $n=2$. Set $E'=\sqrt{3g_2}$. Then $q_1 =q_3=0$,
\begin{align}
& \textstyle c(E)=E^2-\frac{3}{2}g_2, \quad a_0 (E) = 3E.
\end{align}
(iii) The case $n=3$. Set $E'=0$. Then $q_1 =q_3=0$,
\begin{align}
& \textstyle c(E)=E^3-\frac{45}{4}g_2 E-\frac{135}{4}g_3, \quad a_0 (E) = 6(E^2-\frac{15}{4}g_2).
\end{align}
\end{exa}

\section{Bethe Ansatz and Hermite-Krichever Ansatz} \label{sec:HKA}

We review validity of Bethe Ansatz and Hermite-Krichever Ansatz for the Lam\'e equation, and provide examples of the Hermite-Krichever Ansatz.

The method ``Bethe Ansatz'' appears frequently in physics.
In this manuscript, we use ``Bethe Ansatz'' in a somewhat restricted sense.
We assume that the eigenfucntion of the Lam\'e operator has a special form as Eq.(\ref{BA}).
Then the Bethe Ansatz method replaces the problem of finding eigenstates and eigenvalues of the Hamiltonian with a problem of solving transcendental equations for a finite number of variables which are called the Bethe Ansatz equations (see Eq.(\ref{BAeq})).
We review a proposition that almost all eigenfunction is expressed as the form of Bethe Ansatz. For the case of Lam\'e equation it was discussed in \cite{WW}, and extended to the case of the Heun equation in \cite{GW,Tak1}. Note that it is applied to show validity of the Hermite-Krichever Ansatz.
\begin{prop}
(i) For the case $Q(E)\neq 0$, there exists $t_1, \dots t_n$ and $A$ such that $t_j\neq t_{j'}$ $(j\neq j')$, $t_j \not \in \frac{1}{2} \Zint + \frac{\tau}{2} \Zint $ and $\Lambda (x,E)$ is expressed as
\begin{align}
& \Lambda (x,E)= A \frac{\prod_{j=1}^n \sigma(x+t_j)}{\sigma(x)^{n}\prod_{j=1}^n \sigma(t_j)}\exp \left(-x\sum_{i=1}^n \zeta(t_j)\right), \label{BA} 
\end{align}
where $\sigma (x)$ is the Weierstrass sigma function.\\
(ii) 
The function 
\begin{equation}
\tilde{\Lambda }(x)=\frac{\prod_{j=1}^n \sigma(x+t_j)}{\sigma(x)^{n}\prod_{j=1}^n \sigma(t_j)}\exp(cx),
\label{BV}
\end{equation}
with the condition $t_j\neq t_{j'}$ $(j\neq j')$ and $t_j \not \in \frac{1}{2} \Zint + \frac{\tau}{2} \Zint $ is an eigenfunction of the Lam\'e operator (see Eq.(\ref{LameH})), if and only if  $t_j$ $(j=1,\dots ,n)$ and $c$ satisfy the relations,
\begin{align}
& -n\zeta(t_j) + \sum_{k\neq j} \zeta(t_j - t_k)= c\; \; (j=1, \dots ,n).
\label{BAeq} 
\end{align}
The eigenvalue $E$ is given by
\begin{align}
& E=-c^2-n\sum_{j=1}^n(\wp(t_j)-\zeta(t_j)^2)+\sum_{j<k}(\wp(t_j-t_k)-\zeta(t_j-t_k)^2).  \label{BAE} 
\end{align}
\end{prop}

Now we deal with the Hermite-Krichever Ansatz. In our situation, the Hermite-Krichever Ansatz asserts that the differential equation has solutions that are expressed as a finite series in the derivatives of an elliptic Baker-Akhiezer function, multiplied by an exponential function. More precisely, we are going to find solutions to Eq.(\ref{eqn:lame}) of the form
\begin{align}
& f(x) = \exp \left( \kappa x \right) \left( \sum_{j=0}^{n-1} \tilde{b} _j \left( \frac{d}{dx} \right) ^{j} \Phi (x, \alpha ) \right)
\label{Lalpha}
\end{align}
where $\Phi (x,\alpha )= \exp (\zeta( \alpha )x) \sigma (\alpha -x) /(\sigma (x) \sigma (\alpha ))$.
From Eq.(\ref{Lalpha}) we have 
\begin{align} 
& f(x+1,E)  = \exp (-2\eta _1 \alpha +\zeta (\alpha ) +\kappa ) \label{ellint} ,\\
& f(x+\tau,E)  = \exp (-2\eta _3 \alpha +\tau( \zeta (\alpha ) + \kappa ) ) \label{ellintt}.
\end{align}
If the value $\wp (\alpha )$ is calculated, then the values $\alpha $ and $\zeta (\alpha )$ are expressed as elliptic integrals. Hence, if a solution $f(x)$ is expressed as Eq.(\ref{Lalpha}), and $\wp (\alpha )$ and $\kappa $ are calculated explicitly, then we obtain a monodromy formula by elliptic integrals.

Hermite and Halphen investigated the Lam\'e equation for the cases $n=1,2,3,4,5$ on this Ansatz in the 19th century. Belokolos, Eilbeck, Enolskii, Kostov and Smirnov studied for the case $l_0=2, l_1=1, l_2=0, l_3=0$ and the case $l_0=2, l_1 =1 ,l_2 =1, l_3=0$ on the Heun equation (see \cite{BE} and the reference therein). 
On the other hand, Treibich and Verdier \cite{TV} constructed a theory of elliptic soliton following Krichever's idea and studied ``tangential covering'' which would be closely related with the Hermite-Krichever Ansatz.

Now we remind validity of the Hermite-Krichever Ansatz.
\begin{thm} \label{thm:alpha} \cite{Tak4}
Assume $n\in \Zint_{\geq 1}$. There exists polynomials $P_1 (E) ,\dots ,P_4(E)$ such that, if $P_2(E) \neq 0$, then a solution to Eq.(\ref{eqn:lame}) is written as Eq.(\ref{Lalpha}) for some values $\alpha $, $\kappa $ and $\tilde{b} _j$ $(j= 0,\dots ,n-1)$. The values $\alpha $ and $\kappa $ are expressed as
\begin{equation}
 \wp (\alpha ) =\frac{P_1 (E)}{P_2 (E)}, \; \; \; \kappa  =\frac{P_3 (E)}{P_4 (E)} \sqrt{-Q(E)}.
\label{P1P6}
\end{equation}
\end{thm}
Note that $\alpha $ in Eq.(\ref{Lalpha}) and $t_j$ in Eq.(\ref{BA}) satisfy the relation $\alpha = -\sum _{j=1}^n t_j$, which is crucial to show validity of the Hermite-Krichever Ansatz, and the function $\Lambda (x,E)$ in Eq.(\ref{eqn:Lam}) is also expressed as the r.h.s. of Eq.(\ref{Lalpha}).

To calculate the polynomials $P_1 (E) ,\dots ,P_4(E)$, it is effective to apply notions ``twisted Lam\'e polynomial'' and ``theta-twisted Lam\'e polynomial'' which were invented by Maier \cite{Mai}. They are generalized to the case of the Heun equation in \cite{Tak4}.

\begin{exa}
(i) The case $n=1$. 
\begin{align}
& \wp( \alpha )= -E , \quad \kappa =0.
\label{wpn1}
\end{align}
(ii) The case $n=2$. 
\begin{align}
& \wp( \alpha )= e_1 -\frac{(E-3e_1)(E+6e_1)^2}{9(E^2-3g_2)}, \quad \kappa =\frac{2}{3(E^2-3g_2)}\sqrt{-Q(E)}.
\end{align}
Note that the term $(E+6e_1)$ corresponds to the twisted Lam\'e polynomial.\\
(iii) The case $n=3$.
\begin{align}
& \wp( \alpha )= e_1 -\frac{(E^2+6e_1E+45e_1^2-15g_2)(E^2+15e_1E-225e_1^2+\frac{75}{4}g_2)^2}{36E(E^2-\frac{75}{4}g_2)^2},\\
& \kappa =\frac{5}{6E(E^2-\frac{75}{4}g_2)}\sqrt{-Q(E)}.
\end{align}
Note that the terms $(E^2+15e_1E-225e_1^2+\frac{75}{4}g_2)$ and $(E^2-\frac{75}{4}g_2)$ correspond to the twisted Lam\'e polynomials.
\end{exa}

\section{Hyperelliptic-to-elliptic reduction integral formula} \label{sec:intf}

We obtained two expressions of the monodromy in Eqs.(\ref{hypellint}, \ref{hypellinttau}) and Eqs.(\ref{ellint}, \ref{ellintt}). By comparing these two expressions we obtain hyperelliptic-to-elliptic reduction integral formulae. For the proof, see \cite{Tak4}.
\begin{thm} \cite[\S 3]{Tak4}
Let $P_1 (E) , \dots ,P_4 (E)$ be the polynomials defined in Theorem \ref{thm:alpha} and let $a_0(E), c(E)$ be the ones in section \ref{sec:prod}. Set
\begin{align}
& \xi  = \frac{P_1 (E)}{P_2(E)}, \quad \kappa = \frac{P_3(E)}{P_4(E)} \sqrt{-Q(E)},\label{wpas4} 
\end{align}
(i) We have a formula
\begin{equation}
-\frac{1}{2} \int _{\infty}^{E} \frac{a(\tilde{E})}{\sqrt{-Q(\tilde{E})}}d\tilde{E} = \int _{\infty} ^{\xi } \frac{d\tilde{\xi } }{\sqrt{4\tilde{\xi } ^3-g_2\tilde{\xi } -g_3}},
\label{alpints4}
\end{equation}
which reduces a hyperelliptic integral of the first kind to an elliptic integral of the first kind.\\
(ii) Let $E'$ be the value that satisfies $Q(E')=0$ and $P_2(E') \neq 0$, and let $\xi '$ be the value $\xi $ (see Eq.(\ref{wpas4})) evaluated at $E=E'$. Then we have a formula
\begin{equation}
\frac{1}{2} \int  _{E'}^{E} \frac{c(\tilde{E})}{\sqrt{-Q(\tilde{E})}}d\tilde{E} = -\kappa + \int _{\xi ' } ^{\xi } \frac{\tilde{\xi } d\tilde{\xi } }{\sqrt{4\tilde{\xi } ^3-g_2\tilde{\xi } -g_3}} , \label{kap123s4}
\end{equation}
which reduces a hyperelliptic integral of the second kind to an elliptic integral of the second kind.
\end{thm}

\begin{exa}
(i) The case $n=1$. The values $\xi $ and $\kappa $ are given by $\xi = -E $ and $\kappa =0$. Hence the formulae (\ref{alpints4}, \ref{kap123s4}) are trivial.\\
(ii) The case $n=2$. The values $\xi $ and $\kappa $ are given by 
\begin{align}
& \xi = e_1 -\frac{(E-3e_1)(E+6e_1)^2}{9(E^2-3g_2)}, \quad \kappa =\frac{2}{3(E^2-3g_2)}\sqrt{-(E^2-3g_2)\prod _{i=1}^3 (E-3e_i)}.
\end{align}
Set $E' =3e_1$. Then we have $\xi '=e_1$. 
Eqs.(\ref{alpints4}, \ref{kap123s4}) are written as
\begin{align}
& -\frac{1}{2} \int _{\infty}^{E} \frac{3\tilde{E}}{\sqrt{-(\tilde{E}^2 -3g_2)\prod _{i=1}^3 (\tilde{E}-3e_i)}}d\tilde{E} = \int _{\infty} ^{\xi } \frac{d\tilde{\xi } }{\sqrt{4\tilde{\xi } ^3-g_2\tilde{\xi } -g_3}}, \label{eqn201}\\
& \frac{1}{2} \int  _{3e_1}^{E} \frac{\tilde{E} ^2-\frac{3}{2} g_2}{\sqrt{-(\tilde{E}^2 -3g_2)\prod _{i=1}^3 (\tilde{E}-3e_i)}}d\tilde{E} = -\kappa + \int _{e_1 } ^{\xi } \frac{\tilde{\xi } d\tilde{\xi } }{\sqrt{4\tilde{\xi } ^3-g_2\tilde{\xi } -g_3}} \label{eqn202}.
\end{align}
These formulae reduce hyperelliptic integrals of genus two to elliptic integrals.
By setting $\xi =-y/6$, $E=z$, $g_2=a/3$, $g_3=b/54$, we recover the formula 
\begin{equation}
\int \frac{zdz}{\sqrt{(z^2-a)(8z^3-6az-b)}} = \frac{1}{2\sqrt{3}} \int \frac{dy}{\sqrt{y^3-3ay+b}}, \quad \left( y=\frac{2z^3-b}{3(z^2-a)} \right) , \label{Hermitef}
\end{equation}
from Eq.(\ref{eqn201}). This formula was found by Hermite in the 19th century. From Eq.(\ref{eqn202}) we have
\begin{equation}
\int \frac{(2z^2-a)dz}{\sqrt{(z^2-a)(8z^3-6az-b)}} -\frac{1}{3}\sqrt{\frac{8z^3-6az-b}{z^2-a}} = \frac{1}{2\sqrt{3}} \int \frac{ydy}{\sqrt{y^3-3ay+b}}. \label{Hermites}
\end{equation}
(iii) The case $n=3$. The values $\xi $ and $\kappa $ are given by 
\begin{align}
& \xi = e_1 -\frac{(E^2+6e_1E+45e_1^2-15g_2)(E^2+15e_1E-225e_1^2+\frac{75}{4}g_2)^2}{36E(E^2-\frac{75}{4}g_2)^2},\\
& \kappa =\frac{5}{6E(E^2-\frac{75}{4}g_2)}\sqrt{- E\prod _{i=1}^3 (E^2+6e_i E+45e_i ^2-15g_2)}.
\end{align}
Let $E'$ be a root of an equation $E^2+6e_1E+45e_1^2-15g_2=0$. Then we have $\xi '=e_1$. Eqs.(\ref{alpints4}, \ref{kap123s4}) are written as
\begin{align}
& -\frac{1}{2} \int _{\infty}^{E} \frac{6(\tilde{E}^2-\frac{15}{4}g_2)}{\sqrt{-\tilde{E}\prod _{i=1}^3 (\tilde{E}^2+6e_i \tilde{E}+45e_i ^2-15g_2)}}d\tilde{E} = \int _{\infty} ^{\xi } \frac{d\tilde{\xi } }{\sqrt{4\tilde{\xi } ^3-g_2\tilde{\xi } -g_3}},\\
& \frac{1}{2} \int  _{E'}^{E} \frac{\tilde{E} ^3-\frac{45}{4} g_2 \tilde{E} -\frac{135}{4} g_3}{\sqrt{-\tilde{E}\prod _{i=1}^3 (\tilde{E}^2+6e_i \tilde{E}+45e_i ^2-15g_2)}}d\tilde{E} = -\kappa + \int _{e_1 } ^{\xi } \frac{\tilde{\xi } d\tilde{\xi } }{\sqrt{4\tilde{\xi } ^3-g_2\tilde{\xi } -g_3}} .
\end{align}
These formulae reduce hyperelliptic integrals of genus three to elliptic integrals.
\end{exa}

\section{Boundary value problems and Finite-gap potential} \label{sec:bound}

We consider boundary value problems of the Lam\'e operator $H$.
Let $\sigma _{int} (H)$ be the set of eigenvalues of $H$ whose eigenvector is square-integrable on the interval $(0,1)$, i.e.
\begin{equation}
\sigma _{int} (H) =\{ E \: | \: \exists f(x) \neq 0 \mbox{ s.t. }  Hf(x)=Ef(x), \int _0^1 |f(x)| dx < + \infty \}.
\end{equation}
Let $\sigma _d (H)$ be the set of eigenvalues of $H$ whose eigenvector is doubly-periodic, i.e.
\begin{align}
& \sigma _d (H) =\{ E \: | \: \exists f(x) \neq 0 \mbox{ s.t. }  Hf(x)=Ef(x), \; f(x+1)=\pm f(x) ,  \; f(x+\tau )=\pm f(x) \}, 
\end{align}
Note that the doubly-periodic eigenvector is simply the Lam\'e polynomial up to gauge transformation and changing variable. It is known that $\# \sigma _d (H) =2n+1$ and 
\begin{align}
& \sigma _d (H) = \{ E \: | \: Q(E)=0  \}.
\end{align}
Let $\sigma _s (H)$ be the set of eigenvalues of $H$ whose eigenvector is singly-periodic, i.e. 
\begin{align}
& \sigma _s (H) =\{ E \: | \: \exists f(x) \neq 0 \mbox{ s.t. }  Hf(x)=Ef(x), f(x+1)=\pm f(x) \}, 
\end{align}
On the sets $\sigma _{int} (H)$, $\sigma _d (H)$ and $\sigma _s (H)$, we have

\begin{prop} \cite{Taka} \label{prop:spec1}
Assume that $\tau \in \sqrt {-1} \Rea _{>0}$. Then
\begin{eqnarray}
& \sigma _{int} (H) \coprod \sigma _d (H) = \sigma _s (H) .&
\end{eqnarray}
\end{prop}

Next, we briefly explain a relationship with finite-gap potential.
Let 
\begin{align}
& H_s = -\frac{d ^2}{d x^2} +n(n+1) \wp (x+\tau /2) 
\label{eq:lameI}
\end{align}
be the operator which is obtained by the shift $x \rightarrow x+\tau/2$. Then the potential of the operator $H_s$ does not have poles on $\Rea$ and it is real-analytic, if $\tau \in \sqrt {-1} \Rea _{>0}$. Let $\sigma _b(H_s)$ be the set such that 
$$
E \in \sigma _b(H_s) \; \Leftrightarrow \mbox{ All solutions to }(H_s-E)f(x)=0 \mbox{ are bounded on }x \in \Rea .
$$
It is known in the theory of Hill's equation \cite{MW} that, if the potential of the operator $-d^2/dx^2 +q(x) $ is real-analytic and singly-periodic, then the open set $\Rea \setminus \overline{\sigma _b(-d^2/dx^2 +q(x) )}$ is expressed as $(-\infty,E_{0}) \cup \cup _{i\in I} (E_{2i-1},E_{2i})$ where $I= \{ 1, \dots ,n\}$ $(n \in \Zint_{\geq 0})$ or $I= \Zint$. If the number of the unstable bands (gaps) is finite (i.e. $I= \{ 1, \dots ,n\}$ for some $(n \in \Zint_{\geq 0})$), then the potential $q(x)$ is called finite-gap potential or finite-band potential.
 
It is known that the potentials of the Lam\'e operator (see Eq.(\ref{eq:lameI})) for $n\in \Zint_{\geq 1}$ is typical examples of the finite-gap potential, which was established by Ince \cite{I}.
\begin{thm}
Assume that $\tau \in \sqrt {-1} \Rea _{>0}$ and $n \in \Zint_{\geq 1}$. Then  
\begin{equation}
\Rea \setminus \overline{\sigma _b(H_s)}= (-\infty,E_{0}) \cup (E_{1},E_{2}) \cup \dots \cup (E_{2n-1}, E_{2n})
\end{equation}
where $E_0<E_{1}<\cdots <E_{2n}$ and $E_i \in \sigma _d (H)$ $(i=0,\dots ,2n)$.
\end{thm}
Now we present a brief sketch of a new proof of the theorem which is based on the monodromy formule in terms of hyperelliptic integral.
From Eq.(\ref{hypellint}) we have $\Lambda (x+1,E) = \pm \exp (-\tilde{m}_1/2 )  \Lambda (x,E) $, where 
\begin{equation}
\tilde{m}_1 = \int_{E'}^{E} \frac{c(\tilde{E})  -2\eta _1 a_0 (\tilde{E}) }{\sqrt{-Q(\tilde{E})}} d\tilde{E},
\end{equation}
and $E'$ satisfies $Q(E')=0$. Note that, if $\tau \in \sqrt {-1} \Rea _{>0}$, then $a_0(E), c(E), Q(E) \in \Rea$ for $E \in \Rea$. If $Q(E)>0$, then the integrand is pure imaginary by choosing $E'$ appropriately, and we have an equality $|\exp (-\tilde{m}_1/2) |=1$. Thus the eigenfunction is bounded. If $Q(E)<0$, then the integrand is real by choosing $E'$ appropriately, and unboundedness of the eigenfunction is shown.

In the theory of stationary soliton equation (see \cite{GH} and the reference therein), alternative definition of the finite-gap potential is known.
If there exists an odd-order differential operator 
$A= \left( d/dx \right)^{2g+1} +  \! $ $ \sum_{j=0}^{2g-1}\! $ $ b_j(x) \left( d/dx \right)^{2g-1-j} $ such that $[A, -d^2/dx^2+q(x)]=0$, then $q(x)$ is called the algebro-geometric finite-gap potential. Note that the equation  $[A, -d^2/dx^2+q(x)]=0$ is equivalent to the function $q(x)$ being a solution to some stationary higher-order KdV equation.
It is known that, if $q(x) $ is real-analytic on $\Rea$ and $q(x+1)=q(x)$, then $q(x)$ is a finite-gap potential if and only if $q(x)$ is an algebro-geometric finite-gap potential.
For the case of the Lam\'e equation (and the Heun equation), the operator $A$ is calculated by using the function $\Xi (x,E)$.
Write
\begin{equation}
\Xi(x,E) = \sum_{i=0}^{n} \tilde{a}_{n-i}(x) E^i. \label{Xiag}
\end{equation}
From Proposition \ref{prop:prod} we have $\tilde{a}_0(x)=1$.
\begin{prop} \cite{Tak3} \label{thm:PhiA}
Define the $(2n+1)$st-order differential operator $A$ by
\begin{equation}
A= \sum_{j=0}^{n} \left\{ \tilde{a}_j(x+\tau/2)\frac{d}{dx}-\frac{1}{2} \left( \frac{d}{dx} \tilde{a}_j(x+\tau/2) \right) \right\} H_s ^{n-j}. \label{Adef}
\end{equation}
Then the operator $A$ commutes with the operator $H_s$.
Moreover we have
\begin{equation} 
A^2+Q(H_s)=0. \label{algrel}
\end{equation}
\end{prop}
Another formula for the operator $A$ of determinant type is obtained in \cite[\S 3.2]{Tak3}.

\section{The case $n=1$ and analytic continuation of eigenvalues} \label{sec:n1}

In this section we consider the case $n=1$. The Hamiltonian is given as
\begin{equation}
H = -\frac{d ^2}{d x^2} +2\wp (x) . \label{LameHn1}
\end{equation}
Eigenvalues which admits doubly-periodic eigenfunction are written as
\begin{equation}
\sigma _d (H) =\{-e_1, -e_2, -e_3 \}.
\end{equation}
The functions $e_1$, $e_2$ and $e_3$ are analytic on the upper half plane $\Rea + \sqrt{-1} \Rea _{>0}$ and they are related with automorphic functions.

Next we consider eigenvalues which admits square-integrable eigenfunctions.
Set $p=\exp(\pi \sqrt{-1} \tau)$. Then $\tau \rightarrow \sqrt{-1} \infty$ corresponds to $p\rightarrow 0$, and $\tau \in \Rea + \sqrt{-1} \Rea _{>0}$ corresponds to $|p|< 1$.
For the case $p=0$, a complete orthonormal basis for square-integrable eigenfunctions are written as $v_m$ $(m \in \Zint _{\geq 0})$ in Eq.(\ref{Gegenpol}) with the eigenvalue $E_m =\pi ^2 (-2/3+(m+2)^2)$.
Based on the eigenvalues $E_m$ and the eigenfunctions $v_m$ for the case $p=0$, we determine eigenvalues $E_m (p) = E_m+\sum_{k=1}^{\infty} E_m^{\{2k\}}p^{2k}$ and normalized eigenfunctions $v_m(p)= v_m+ \sum_{k=1}^{\infty} \sum_{m'\in \Zint _{\geq 0}} c_{m ,m'}^{\{ 2k \} }v_{m'}p^{2k}$ for the operator $H$
as formal power series in $p$ by applying a method of perturbation.
For details see \cite{Tak2,Taka}.
Convergence of the formal power series of eigenvalues in the variable $p$ was shown in \cite{Tak2}.
\begin{prop} \label{cor:conv} \cite[Corollary 3.7]{Tak2}
Let $E_m(p)$ $(m\in \Zint _{\geq 0})$ (resp. $v_m(p)$) be the formal eigenvalue (resp. eigenfunction) of the Hamiltonian $H$.
If $|p|$ is sufficiently small then the power series $E_m(p)$ converges.
\end{prop}
Note that it was also shown in \cite{Tak2} that the eigenvalues $E_m (p)$ are real-analytic in $p^2 \in (-1,1)$.

The inferred convergence radius and expansions of the first few terms for the eigenvalue $E_m(p)$ are calculated as follows (see also \cite{Taka}):

\begin{center}
\begin{tabular}{|l|l|l|} \hline
$E_0 (p)$ & $\pi^2 \left( \frac{10}{3} + \frac{80}{3} p^2+\frac{1360}{27} p^4+\frac{20800}{243} p^6 +\frac{195920}{2187}p^8+ \frac{3174880}{19683} p^{10} + \dots \right)$ & .749 \\ \hline
$E_2 (p)$ & $\pi^2 \left( \frac{46}{3} + \frac{272}{15} p^2+\frac{198928}{3375} p^4+\frac{55403584}{759375} p^6 +\frac{4307155408}{34171875}p^8 +\dots \right)$ & .749 \\ \hline
$E_4 (p)$ & $\pi^2 \left( \frac{106}{3} + \frac{592}{35} p^2 +\frac{2279248}{42875} p^4+\frac{3773733184}{52521875} p^6 +\frac{1634762851088}{12867859375}p^8 +\dots \right)$ & .875 \\ \hline \hline
$E_1 (p)$ & $\pi^2 \left( \frac{25}{3} + 20 p^2+ 65 p^4+\frac{115}{2} p^6 +\frac{2165}{16}p^8+ \frac{3165}{32} p^{10} +\frac{23965}{128}p^{12} +\dots \right)$ & .838 \\ \hline
$E_3 (p)$ & $\pi^2 \left( \frac{73}{3} + \frac{52}{3} p^2 +\frac{1493}{27} p^4+\frac{35671}{486} p^6 +\frac{4492153}{34992}p^8 + \frac{55853449}{629856} p^{10} + \dots \right)$ & .838 \\ \hline
$E_5 (p)$ & $\pi^2 \left( \frac{241}{3} + \frac{82}{5} p^2+\frac{50339}{1000} p^4+\frac{13640101}{200000} p^6 +\frac{3872868499}{32000000}p^8 +\dots \right)$ & .906 \\ \hline
\end{tabular} 
\\
\vspace{.05in}
\noindent {\it Table 1. Expansion and inferred convergence radius.} 
\end{center}

If the function $E_0 (p)$ is analytic on the upper half plane $\Rea + \sqrt{-1} \Rea _{>0}$, then the radius of convergence in $p$ should be $1$, thought it is numerically inferred to be $.749$. Hence it seems that the function $E_0 (p)$ has singularity.

Now we consider analytic continuation of the eigenvalues $E_m (p)$ $(m \in \Zint _{\geq 0})$ in the variable $p$. From Eq.(\ref{hypellint}) or Eqs.(\ref{ellint}, \ref{wpn1}),  the eigenvalue $E$ is continued analytically in $p$ 
by keeping the conditions 
\begin{align}
& E = -\wp(t_0), \label{acc23} \\
& \exists m \in \Zint, \; \; \;  -\zeta (t_0) +2 \eta _1 t_0 = m \pi \sqrt{-1}. \label{acc24}
\end{align}
Let $\Re a $ (resp. $\Im a$) be the real part (resp. the imaginary part) of the number $a$, and let ${\mathcal C}_a$ be the cycle starting from $\Re a$, approaching the point $a$ parallel to the imaginary axis, turning anti-clockwise around $a$ and returning to $\Re a$ as shown in Figure 2.

\begin{center}
\begin{picture}(200,100)(0,0)
\put(10,10){\vector(1,0){170}}
\put(10,10){\vector(0,1){85}}
\put(112,10){\vector(0,1){40}}
\put(112,50){\line (0,1){19}}
\put(108,69){\vector(0,-1){40}}
\put(108,30){\line(0,-1){20}}
\put(110,75){\circle*{2}}
\put(108,75){\oval(10,12)[l]}
\put(112,75){\oval(10,12)[r]}
\put(112,81){\vector(-1,0){4}}
\put(185,6){Re}
\put(15,85){Im}
\put(108,1){$\Re a$}
\put(118,68){$a$}
\put(96,45){${\mathcal C}_a$}
\end{picture}
\\
\vspace{.05in}
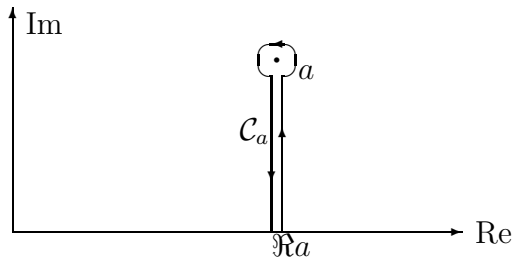
\noindent {\it Figure 2. Cycle ${\mathcal C}_a$.}
\end{center}

We continue the eigenvalue $E$ analytically along the cycle ${\mathcal C}_a$ on where $a$ is a possible branching point for $|a|<.90$, $\Re a \geq 0$ and $\Im a \geq 0$. The branching points for $|a|<.90$, $\Re a \geq 0$ and $\Im a \geq 0$ are calculated in \cite{Taka}, and they are listed in Table 3. Branching along the cycle ${\mathcal C}_a$ is determined as shown in Table 3 (see also \cite{Taka}).

\begin{center}
\begin{small}
\begin{tabular}{|l|l|}
\hline
$a=.258666+.697448 I$ & $E_0 (p)  \Rightarrow E_2 (p), \; E_2 (p) \Rightarrow E_0 (p), \: \; E_4 (p) \Rightarrow E_4 (p), \; E_6 (p) \Rightarrow E_6 (p)$ \\
$a=.224582 +.842777 I$ & $E_0 (p) \Rightarrow E_4 (p), \; E_2 (p) \Rightarrow E_2 (p), \: \; E_4 (p) \Rightarrow E_0 (p), \; E_6 (p) \Rightarrow E_6 (p)$ \\
$a=.552288 +.677536 I$ & $E_0 (p) \Rightarrow E_4 (p), \; E_2 (p) \Rightarrow E_2 (p), \: \; E_4 (p) \Rightarrow E_0 (p), \; E_6 (p) \Rightarrow E_6 (p)$ \\
$a=.314813+.821858 I$ & $E_0 (p) \Rightarrow E_4 (p), \; E_2 (p) \Rightarrow E_2 (p), \: \; E_4 (p) \Rightarrow E_0 (p), \; E_6 (p) \Rightarrow E_6 (p)$ \\
$a=.686317+.559106 I$ & $E_0 (p) \Rightarrow E_0 (p), \; E_2 (p) \Rightarrow E_4 (p), \: \; E_4 (p) \Rightarrow E_2 (p), \; E_6 (p) \Rightarrow E_6 (p)$ \\
\hline
$a=.535905+.640487 I$ & $E_1 (p) \Rightarrow E_3(p) , \; E_3(p)  \Rightarrow E_1(p) , \: \; E_5(p)  \Rightarrow E_5(p) , \; E_7(p)  \Rightarrow E_7(p) $ \\
\hline
\end{tabular}
\end{small}
\vspace{.05in}
\\
\noindent {\it Table 3. Branching along the cycle ${\mathcal C}_a$} 
\end{center}

The closest branching point from the origin is $p=.258666+.697448 I$ ($|p|=.743869$) and the eigenvalues $E_0(p)$ and $E_2(p)$ are connected by continuing analytically along the cycle ${\mathcal C}_p$ $(p=.258666+.697448I)$ (see Table 3). 
On the other hand, by the method of perturbation the convergence radii of the expansions of the eigenvalues $E_0(p)$ and $E_2(p)$ are both inferred to be around $.749$ from Table 1.
Thus, convergence radii calculated by different methods are very close and compatibility between the method of perturbation and the method of monodromy is confirmed. Moreover, we obtain a reason why the convergence radii of the eigenvalues $E_0(p)$ and $E_2(p)$ are very close by considering the branching point. For details see \cite{Taka}.

It would be able to consider analytic continuation of eigenvalues for the case $n \geq 2$ and the case of the Heun equation by applying results on the Hermite-Krichever Ansatz.

\end{document}